\documentclass[8pt,a4paper]{article}
\usepackage{amsfonts}
\usepackage{amsmath}
\usepackage{amssymb}
\usepackage{amsthm}
\usepackage[dvips]{graphicx}
\usepackage{latexsym}
\usepackage{mathrsfs}
\usepackage{amscd}
\usepackage{indentfirst}
\usepackage{textcomp}
\usepackage[numbers,square,sort&compress]{natbib} 
\usepackage[colorlinks=true,allcolors=black ]{hyperref}
\allowdisplaybreaks[4]

\newtheorem{Theorem}{Theorem}[section]

\newtheorem{Corollary}[Theorem]{Corollary}
\newtheorem{Lemma}[Theorem]{Lemma}

\newtheorem{Remark}[Theorem]{Remark}
\theoremstyle{definition}
\newtheorem{Definition}{Definition}[section]
\numberwithin{equation}{section}

\DeclareMathOperator{\Der}{Der}
\DeclareMathOperator{\BDer}{BDer}
\DeclareMathOperator{\IBDer}{IBDer}

\DeclareMathOperator{\Span}{span}
\renewcommand{\d}{\mathrm{d}}
\newcommand{\BF}{\mathbb{F}}

\newcommand{\BN}{\mathbb{N}}
\newcommand{\BZ}{\mathbb{Z}}

\begin{document}
\title{{\bf  Skew-symmetric super-biderivations of the special Lie superalgebra $S(m,n;\underline{t})$}}
\author{\normalsize \bf Da  Xu$^1$\,\,\,\,Xiaoning  Xu$^1$}
\date{{{\small{  1. School of Mathematics, Liaoning University, Shenyang, 110036,
China  }}}}
\maketitle

\begin{abstract}
This paper aims to study the skew-symmetric super-biderivations of the special Lie superalgebra $S(m,n;\underline{t})$. Let $S$ denote the special Lie superalgebra $S(m,n;t)$ over a field of characteristic $p>2$. Utilizing the abelian subalgebra $T_S$ and the weight space decomposition of $S$ with respect to $T_S$, we show the action of a skew-symmetric super-biderivation on the elements of $T_S$ and some specific elements of $S$. Moreover, we prove that all skew-symmetric super-biderivations of $S$ are inner.
\end{abstract}
\renewcommand{\thefootnote}{\fnsymbol{footnote}}
\footnote[0]{ Project Supported by National Natural Science
Foundation of China (No.11501274) and the Science Research Project of Liaoning Provincial Education Department, China (No.L2015203).
\\ Author Email: lnuxxn@163.com (X. Xu)}

\section{Introduction}

The research of derivations as pivotal subjects in the investigation of rings and algebras has been underway for many years. Biderivations, as a generation of derivations, have aroused the interest of many scholars in recent years. The concept of a symmetric bi-derivation was first introduced by  Maksa in the real linear space \cite{Maksa1987}. In \cite{Bresar1993}, it was shown that every biderivation $D$ of a noncommutative prime ring $R$ is of the form $D(x,y)=\lambda[x,y]$ for some $\lambda\in \mathbb{C}$. In \cite{zhang2006}, the authors proved that every biderivation of the associated nest algebra $\tau(\mathcal{N})$ is an inner biderivation if and only if $\text{dim}\,0_+\neq 1$ and $\text{dim}\,H^{\perp}_{-}\neq 1$. The notation of biderivations of Lie algebras was introduced and studied in \cite{wang2011}.  Since then biderivations of Lie algebras have been studied widely \cite{chang20191,chang20192,chen2016,cheng2017,han2016,tang20181,tang20182,wang2013}. Scholars proved that each anti-symmetric biderivation is inner for restricted Cartan-type Lie algebras, simple generalized Witt algebras and the Lie algebra $\mathfrak{gca}$ (see \cite{chang20191,chang20192,chen2016,cheng2017}, for details). 

Lie superalgebras as a generalization of Lie algebras came from supersymmetry in mathematical physics. Super-biderivations of Lie superalgebras have naturally aroused the interest of numerous scholars \cite{bai2023,chang2021,cheng2019,Dilxat2023,fan2017,li2018,tang2020,xu2015,yuan2021,yuan2018,zhao2020}. The definition of super-biderivations of Lie superalgebras was introduced for the first time in the Heisenberg superalgebra \cite{xu2015}. In \cite{fan2017}, the author modified some results about super-biderivations of Lie superalgebras and proved that each super-biderivation on the centerless super-Virasoro algebras is inner super-biderivation. Thereafter, Chang \cite{chang2021} and Zhao \cite{zhao2020} proved that all skew-symmetric super-biderivations of generalized Witt Lie superalgebras and contact Lie superalgebras are inner. In \cite{bai2023}, the authors obtained that any super-skew-symmetric superbiderivation of simple modular Lie superealgebra of Witt type and special type is inner. The difference between the aforementioned article \cite{bai2023} and this paper lies in the basic fields: the former necessitates an algebraically closed field, while the latter does not require such a constraint. 

This paper is devoted to studying the skew-symmetric super-biderivations of the special Lie superalgebra $S(m,n;\underline{t})$. And the paper is arranged as follows. In Section 2, we review the basic definitions concerning the special Lie superalgebra $S(m,n;\underline{t})$. In Section 3, we give the definition of skew-symmetric super-biderivations of Lie superalgebras and obtain some useful conclusions about the skew-symmetric super-biderivations. In Section 4, we use the method of the weight space decomposition of $S(m,n;\underline{t})$ with respect to the abelian subalgebra $T_S$ to prove that all skew-symmetric super-biderivations of $S(m,n;\underline{t})$ are inner.

\section{Preliminaries}

In this section, we recall the basic notation concerning the special Lie superalgebras $S(m,n;\underline{t})$ (see \cite{zhang2005}).

Hereafter $\BF$ is an algebraically closed field of characteristic $p>2$ and $\BZ_2 = \{\bar{0}, \bar{1}\}$ is the additive group of order 2. For a vector superspace $V=V_{\bar{0}}\, \oplus\, V_{\bar{1}}$, we write $\d(x)=\alpha$ for the parity of $x\in V_\alpha,\alpha\in \BZ_2$. If $V=\oplus_{i\in \BZ}V_i$ is a $\BZ$-graded vector space, for $x\in V_i,\,i\in \BZ$, $x$ is a $\BZ$-homogeneous element and its $\BZ$-degree $i$ . Throughout this paper, we should mention that once the symbol $\d(x)$  appears in an expression, it implies that $x$ is a $\BZ_2$-homogeneous  element.

Let $\BN_+$ be the set of positive integers and $\BN$ the set of non-negative integers.
Given $m,n\in\BN_+\backslash\{1\}$. For $\alpha=(\alpha_{1},\alpha_{2},\cdots,\alpha_{m})\in \BN^m$, put $|\alpha|=\sum_{i=1}^{m}\alpha_{i}$. For $\beta=(\beta_{1},\beta_{2},\cdots,\beta_{m})\in \BN^m$, we write $\alpha+\beta=(\alpha_{1}+\beta_{1},\alpha_{2}+\beta_{2},\cdots,\alpha_{m}+\beta_{m})$,  $\binom{\alpha}{\beta}=\prod_{i=1}^m\binom{\alpha_i}{\beta_i}$, $\alpha\leq\beta\Longleftrightarrow\alpha_i\leq\beta_i, i=1,2,\cdots,m$. For $\varepsilon _i:=( \delta _{i1},\delta _{i2},\cdots ,\delta _{im})$, where $\delta_{ij}$ is the Kronecker symbol, we abbreviate $x^{(\varepsilon_{i})}$ to $x_{i}$, $i=1, 2, \cdots, m$.  We call $\mathcal{U}(m)$ a $divided\ power\ algebra$ which denotes the  $\BF$-algebra of power series in the variable $x_{1},x_{2},\cdots,x_{m}$. The following formulas hold in $\mathcal{U}(m)$:
  \[x^{(\alpha)}x^{(\beta)}=\binom{\alpha+\beta}{\alpha}x^{(\alpha+\beta)},\ \forall\ \alpha,\beta\in \BN^m.\]
 
Let $\Lambda(n)$ denote the $Grassmann\ superalgebra$ over $\BF$ in $n$ variables $x_{m+1},x_{m+2},\cdots,x_s$, where $s=m+n$. Denote the tensor product $\mathcal{U}(m)\otimes\Lambda(n)$ by $\Lambda(m,n)$. Then $\Lambda(m,n)$ have  a $\BZ_2$-gradation induced by the trivial $\BZ_2$-gradation of  $\mathcal{U}(m)$ and the natural $\BZ_2$-gradation of $\Lambda(n)$:
  \[\Lambda(m,n)_{\bar{0}}=\mathcal{U}(m)\otimes\Lambda(n)_{\bar{0}},\qquad \Lambda(m,n)_{\bar{1}}=\mathcal{U}(m)\otimes\Lambda(n)_{\bar{1}}.\]
  
Obviously, $\Lambda(m,n)$ is an associative superalgebra. For $g\in\mathcal{U}(m), f\in\Lambda(n)$, we simply write $g\otimes f$ as $gf$. The following formulas hold in $\Lambda(m,n)$:
  \[x_ix_j=-x_jx_i, \qquad i,j=m+1,\cdots,s.\]
  \[x^{(\alpha)}x_j=x_jx^{(\alpha)}, \quad\ \forall\ \alpha\in \BN^m,j=m+1,\cdots,s.\]
For $k = 1,\cdots,n$, set
  \[B_k:=\{\langle i_1,i_2,\cdots,i_k\rangle\mid m+1\leq i_1<i_2<\cdots<i_k\leq s\}\]
and $B(n):=\bigcup_{k=0}^nB_k$, where $B_0=\varnothing$. For $u=\langle i_1,i_2,\cdots,i_k\rangle\in B_k$, set $|u|=k, \{u\}=\{i_1,i_2,\cdots i_k\}$ and $x^u=x_{i_1}x_{i_2}\cdots x_{i_k}$. Specially, let $|\varnothing|=0$, $x^\varnothing=1$. It is obvious that $\{x^{(\alpha)}x^u\mid \alpha\in \BN^m,u\in B(n)\}$ is an $\BF$-basis of $\Lambda(m,n)$.

Let $Y_0=\{1,2,\cdots,m\}$, $Y_1=\{m+1,\cdots,s\}$, and $Y=Y_0\cup Y_1$. Let $D_1,D_2,\cdots,D_s$ be the linear transformations of $\Lambda(m,n)$ such that
  \begin{align*}
  	D_i(x^{(\alpha)}x^u)=
  	\begin{cases}
  		x^{(\alpha-\varepsilon_i)}x^u,\quad&\forall\ i\in Y_0,\\
  		x^{(\alpha)}\partial_i(x^u),\quad &\forall\ i\in Y_1,
  	\end{cases}
  \end{align*}
where $\partial_i$ is the special derivation of $\Lambda(n)$. Then $D_1,D_2,\cdots,D_s$ are superderivations of the superalgebra $\Lambda(m,n)$, and $\d(D_i)=\tau(i)$, where
   \begin{align*}
   	\tau(i)=
   	\begin{cases}
   		\bar{0},\qquad \forall\ i\in Y_0,\\
   		\bar{1},\qquad \forall\ i\in Y_1.
   	\end{cases}
   \end{align*}
 Let
  \[W(m,n):=\bigg\{\sum_{i=1}^sf_iD_i\mid f_i\in \Lambda(m,n),\ \forall\ i\in Y\bigg\}.\]
Then $W(m,n)$ is an infinite-dimensional Lie superalgebra which is contained in $\Der(\Lambda(m,n))$ and the following formula holds:
  \[ [fD_i,gD_j]=fD_i(g)D_j-(-1)^{\d(fD_i)\d(gD_j)}gD_j(f)D_i, \]
for all $f,g\in\Lambda(m,n)$ and $i,j\in Y$.

Fix two $m$-tuples of positive integers $\underline{t}=(t_1,t_2,\cdots,t_m)$ and $\pi=(\pi_1,\pi_2,\cdots,\pi_m)$,
where $\pi_i=p^{t_i}-1$ for all $i\in Y_0$ and $p$ is the characteristic of the basic field $\BF$. Let
  \[\Lambda(m,n;\underline{t}):= \Span_{\BF}\{x^{(\alpha)}x^u\mid \alpha\in A(m,\underline{t}),u\in B(n)\}.\]
where $A(m,\underline{t} ) =\{\alpha=(\alpha _1,\alpha _2,\cdots ,\alpha _m)\in \BN^m\mid 0\leq \alpha _i\leq \pi _i,i\in Y_0 \}$. Then $\Lambda(m,n;\underline{t})$ is a subalgebra of $\Lambda(m,n)$. Let 
  \[W(m,n;\underline{t}):=\bigg\{\sum_{i=1}^{s}f_iD_i\mid f_i\in \Lambda(m,n;\underline{t}),\ \forall\ i\in Y\bigg\}.\]
Then $W(m,n;\underline{t})$ is a finite-dimensional simple Lie superalgebra, which is called the generalized Witt Lie superalgebra. Note that $W(m,n;\underline{t})$ possesses a $\BZ$-graded structure:
  \[W(m,n;\underline{t})=\bigoplus\limits_{r=-1}^{\xi-1}W(m,n;\underline{t})_r,\]
by letting $W(m,n;\underline{t})_r:=\Span_{\BF}\{x^{(\alpha)}x^uD_j\mid|\alpha|+|u|=r+1,j\in Y\}$ and $\xi:=|\pi|+n$.

Hereafter, define a linear mapping $D_{t_1t_2}:\Lambda(m,n)\rightarrow W(m,n)$ such that
  \[D_{t_1t_2}(f)=\sum_{i=1}^{2}f_{t_i}D_{t_i},\]
for all $f\in \Lambda(m,n)$, where $f_{t_1}=-(-1)^{\d(f)(\tau(t_1)+\tau(t_2))}D_{t_2}(f)$, 
$f_{t_2}=(-1)^{\tau(t_1)\tau(t_2)}D_{t_1}(f)$ and $\tau(i)=\d(D_i)$ for all $i\in Y$. Let
  \[S(m,n):= \Span_{\BF}\{D_{t_1t_2}(f)\mid t_1,t_2\in Y,f\in \Lambda(m,n)\}.\]
Then $S(m,n)$ is an infinite-dimensional Lie superalgebra. Let
  \[S(m,n;\underline{t}):= \Span_{\BF}\{D_{ij}(f)\mid i,j\in Y,f\in \Lambda(m,n;\underline{t})\}.\]
Then $S(m,n;\underline{t})$ is a finite-dimensional simple Lie superalgebra, which is called the special Lie superalgebra. Note that $S(m,n;\underline{t})$ possesses a $\BZ$-graded structure:
  \[S(m,n;\underline{t})=\bigoplus\limits_{r=-1}^{\xi-2}S(m,n;\underline{t})_r,\]
by letting $S(m,n;\underline{t})_r:=\Span_{\BF}\{D_{ij}(x^{(\alpha)}x^u)\mid|\alpha|+|u|=r+2, i,j\in Y\}$. For convenience, we use $S$ and $S_r$ denote $S(m,n;\underline{t})$ and its $\BZ$-graded subspace $S(m,n;\underline{t})_r$ respectively.

\section{The notion of super-biderivation}

In this section, we give the definition of skew-symmetric super-biderivations of Lie superalgebras and obtain some useful conclusions about the skew-symmetric super-biderivations.

Let $G$ be a Lie algebra over an algebraically closed field. Recall that a linear map $D : G \rightarrow G$ is a derivation of $G$ if
  \[D([x,y])=[D(x),y]+[x,D(y)],\]
for all $x, y\in G$. And we call a bilinear map $\varphi : G \times G\rightarrow G$ is a biderivation of $G$ if the following axioms are satisfied:
  \begin{align*}
  	\varphi (x,[y,z])&=[\varphi (x,y),z]+[y,\varphi (x,z)],\\
  	\varphi([x,y],z)&=[\varphi(x,z),y]+[x,\varphi(y,z)],
  \end{align*} 
for all $x, y, z\in G$. A biderivation $\varphi$ is called skew-symmetric if $\varphi(x, y) = -\varphi(y, x)$ for all $x, y \in G$.

A Lie superalgebra is a vector superspace $L=L_{\bar{0}}\,\oplus\, L_{\bar{1}}$ with an even bilinear mapping $[\cdot,\cdot]:L\times L\rightarrow L$ satisfying the following axioms:
  \begin{align*}
  	[x,y]&=-(-1)^{\d(x)\d(y)}[y,x],\\
  	[x,[y,z]]&=[[x,y],z]+(-1)^{\d(x)\d(y)}[y,[x,z]],
  \end{align*}
for all $x,y,z\in L$. Recall that a linear map $D:L\rightarrow L$ is a superderivation of $L$ if 
  \[D([x,y])=[D(x),y]+(-1)^{\d(D)\d(x)}[x,D(y)],\]
for all $x,y\in L$. Meanwhile, we write $\Der_{\bar{0}}(L)$ (resp. $\Der_{\bar{1}}(L)$) for the set of all superderivations of $\BZ_2$-degree $\bar{0}$ (resp. $\bar{1}$) of $L$. Let $\Der(L)=\Der_{\bar{0}}(L)\oplus\Der_{\bar{1}}(L).$

A $\BZ_2$-homogeneous bilinear map $\phi$ of $\BZ_2$-degree $\gamma$ of $L$ is a bilinear map such that $\phi(L_\alpha,L_\beta)\subset L_{\alpha+\beta+\gamma}$ for any $\alpha,\beta\in\BZ_2$. Specially, we call $\phi$ is even if $\d(\phi)=\gamma=\bar{0}$.

\begin{Definition}\label{def3.1}
	We call a bilinear map $\phi : L\times L\rightarrow L$ is a skew-symmetric super-biderivation of L if the following axioms are satisfied:
	\begin{align}
		\phi(x,[y,z])&=[\phi(x,y),z]+(-1)^{(\d(\phi)+\d(x))\d(y)}[y,\phi(x,z)], \label{eq3.1}\\
		\phi(x,y)&=-(-1)^{\d(\phi)\d(x)+\d(\phi)\d(y)+\d(x)\d(y)}\phi(y,x),\nonumber
	\end{align}
	for all $x,y,z\in L$. 
\end{Definition}
\noindent Observe that if $\phi$ is a skew-symmetric super-biderivation, the equation below is obviously fulled. 
\begin{align}
	&\phi([x,y],z)=[x,\phi(y,z)]+(-1)^{(\d(\phi)+\d(z))\d(y)}[\phi(x,z),y]. \label{eq3.2}
\end{align}
\noindent Meanwhile, we write $\BDer_{\bar{0}}(L)$ (resp. $\BDer_{\bar{1}}(L)$) for the set of all skew-symmetric super-biderivations of $\BZ_2$-degree $\bar{0}$ (resp. $\bar{1}$) of $L$. Let $\BDer(L)=\BDer_{\bar{0}}(L)\oplus\BDer_{\bar{1}}(L).$

\begin{Lemma}\label{lem3.1}
	If the bilinear map $\phi_{\lambda}$ : $L\times L\rightarrow L$ is defined by
        \[\phi_{\lambda}( x,y ) =\lambda [ x,y],\]
	for all $x,y\in L$, where $\lambda \in \BF$, then $\phi_{\lambda}$ is a skew-symmetric super-biderivation of $L$. This class of super-biderivations is called inner. Write $\IBDer(L)$ for the set of all inner super-biderivations of $L$.
	
	\begin{proof}
		Due to $\phi_{\lambda}( x,y ) =\lambda [ x,y]$, we deduce that $\phi_{\lambda}$ is an even super-biderivation, i.e. $\d(\phi_\lambda)=\bar{0}$. By the skew-symmetry of Lie superalgebras, it is easy to see that 
		  \[\phi_\lambda(x,y)=-(-1)^{\d(\phi_\lambda)\d(x)+\d(\phi_\lambda)\d(y)+\d(x)\d(y)}\phi_\lambda(y,x),\]
		for any $x, y \in L$. And it is easily obtain the following three equalities:
		\begin{align*}
			\phi_\lambda(x,[y,z])&=\lambda[x,[y,z]],\\
			[\phi_\lambda(x,y),z]&=\lambda[[x,y],z],\\
			(-1)^{(\d(\phi_\lambda)+\d(x))\d(y)}[y,\phi_\lambda(x,z)]&=(-1)^{\d(x)\d(y)}\lambda[y,[x,z]].
		\end{align*}
		Since the graded Jacobi identity $[x,[y,z]]=[[x,y],z]+(-1)^{\d(x)\d(y)}[y,[x,z]]$,  we have that
		  \[\phi_\lambda(x,[y,z])=[\phi_\lambda(x,y),z]+(-1)^{(\d(\phi_\lambda)+\d(x))\d(y)}[y,\phi_\lambda(x,z)],\]
		for any $x,y,z \in L$. Similarly, it follows that     
		  \[\phi_\lambda([x,y],z)=[x,\phi_\lambda(y,z)]+(-1)^{(\d(\phi_\lambda)+\d(z))\d(y)}[\phi_\lambda(x,z),y],\] 
	    for any $x, y, z \in L$. The proof is completed.
	\end{proof}
   
\end{Lemma}
 
\begin{Lemma}\label{lem3.2}		
	Let $L$ be a Lie superalgebra. Suppose that $\phi$ is a skew-symmetric super-biderivation of $L$. Then
		\[[\phi(x,y),[u,v]]=(-1)^{\d(\phi)(\d(y))+\d(u)}[[ x,y] ,\phi( u,v)],\]
	for any $x,y,u,v\in L$.

	\begin{proof}
		Due to the Definition \ref{def3.1}, there are two different ways to compute $\phi([x, u], [y, v])$. From the equation \eqref{eq3.1}, it follows that
		\begin{align*}
			\phi ([x,u],[y,v])
			=\ &[\phi([x,u],y), v]+(-1)^{(\d(\phi)+\d(x)+\d(u))\d(y)}[y,\phi([x,u],v)] \\
			=\ &[[x,\phi(u,y)], v]+(-1)^{(\d(\phi)+\d(y))\d(u)}[[\phi(x, y), u], v] \\
			&+(-1)^{(\d(\phi)+\d(x)+\d(u))\d(y)}[ y,[ x,\phi (u, v)]] \\
			&+(-1)^{(\d(\phi)+\d(x)+\d(u))\d(y)+(\d(\phi)+\d(v))\d(u)}[y, [\phi(x, v), u]].
		\end{align*}
		According to the equation \eqref{eq3.2}, one gets		
	    \begin{align*}
	    	\phi([x,u],[y,v])
	    	=\ &[x,\phi(u,[y,v])]+(-1)^{(\d(\phi)+\d(y)+\d(v))\d(u)}[\phi(x,[y,v]),u] \\
	    	=\ &[x,[\phi(u,y),v]]+(-1)^{(\d(\phi)+\d(u))\d(y)}[x,[y,\phi(u,v)]] \\
	    	&+(-1)^{(\d(\phi)+\d(y)+\d(v))\d(u)}[[\phi(x,y),v],u] \\
	    	&+(-1)^{(\d(\phi)+\d(y)+\d(v))\d(u)+(\d(\phi)+\d(x))\d(y)}[[y, \phi(x, v)], u].
	    \end{align*}
		
		Comparing two sides of the above two equations, and using the graded Jacobi identity of Lie superalgebras, we have that
		\begin{equation}\label{eq3.3}
			\begin{aligned}
				&[\phi(x,y),[u,v]]-(-1)^{\d(\phi)(\d(y)+\d(u))}[[x,y],\phi(u,v)]\\
				=\ &(-1)^{\d(y)\d(u)+\d(y)\d(v)+\d(u)\d(v)}([\phi(x,v),[u,y]]-(-1)^{\d(\phi)(\d(v)+\d(u))}[[x,v],\phi(u, y)]).
			\end{aligned}		
		\end{equation}		
		And we set 
		  \[f(x,y;u,v)=[\phi(x,y),[u,v]]-(-1)^{\d(\phi)(\d(y)+\d(u))}[[x,y],\phi(u,v)].\]
		According to the equation \eqref{eq3.3}, it easily seen that 
		  \[f(x,y;u,v)=(-1)^{\d(y)\d(u)+\d(y)\d(v)+\d(u)\d(v)}f(x,v;u,y).\]		
		On the one hand, we have
		\begin{align*}
			f(x,y;u,v)&=-(-1)^{\d(u)\d(v)}f(x,y;v,u)\\
			              &=-(-1)^{\d(u)\d(v)}(-1)^{\d(y)\d(v)+\d(y)\d(u)+\d(v)\d(u)}f(x,u;v,y)\\
			              &= (-1)^{\d(y)\d(u)}f(x,u;y,v).
		\end{align*}
		On the other hand, we also get
		\begin{align*}
			f(x,y;u,v)&=(-1)^{\d(y)\d(u)+\d(y)\d(v)+\d(u)\d(v)}f(x,v;u,y)\\
			              &=-(-1)^{\d(y)\d(u)+\d(y)\d(v)+\d(u)\d(v)}(-1)^{\d(u)\d(y)}f(x,v;y,u)\\
			              &=-(-1)^{\d(u)\d(y)}f(x,u;y,v).
		\end{align*}		
		Hence, it follows at once that 
		  \[f(x,y;u,v)=-f(x,y;u,v).\]		
		Due to the characteristic of $\mathbb{F}$ is different from 2, we obtain that $f(x,y;u,v)=0.$ So we have 
		  \[[\phi(x,y),[u,v]]=(-1)^{\d(\phi)(\d(y)+\d(u))}[[x,y],\phi(u, v)].\]
		The proof is completed.
	\end{proof}

\end{Lemma}

\begin{Lemma}\label{lem3.3}		
	Let $L$ be a Lie superalgebra. Suppose that $\phi$ is a skew-symmetric super-biderivation of $L$. If $\d(x)+\d(y)=\bar{0}$, then
		\[[\phi ( x,y ) ,[ x,y]] =0,\]
	for any $x,y\in L$.

	\begin{proof}
		By Lemma \ref{lem3.2}, we can have
			\[[\phi(x,y),[x,y]]=(-1)^{\d(\phi)(\d(y)+\d(x))}[[x,y],\phi(x, y)].\]		
		In view of $\d(x)+\d(y)=\bar{0}$, it is following that
		\begin{align*}
			[\phi(x,y),[x,y]]&=[[x,y],\phi(x,y)]\\
			                       &=-[\phi(x,y),[x,y]].
		\end{align*}
		Therefore, we obtain $[\phi(x,y),[x, y] ]=0$.
	\end{proof}

\end{Lemma}
 
 \begin{Lemma}\label{lem3.4}   	
 	Let $L$ be a Lie superalgebra. Suppose that $\phi$ is a skew-symmetric super-biderivation of $L$. If $[x,y]=0$, then
 	\[\phi(x,y)\in Z_L([L,L]),\]
 	where $Z_L([L,L])$ is the centralizer of  $[L,L]$. 	

 	\begin{proof}
 		If $[x,y]=0$, then we obtain
 		\begin{align*}
 			[\phi(x,y),[u,v] ]&= (-1)^{\d(\phi)(\d(y)+\d(u))}[[x,y],\phi[u,v]]\\
 			                             &= 0,
 		\end{align*}
 		for any $u, v \in L$. So we get that $\phi(x,y)\in Z_L([ L,L])$.
 	\end{proof}
 		
 \end{Lemma}

\section{Skew-symmetric super-biderivation of $S(m,n;\underline{t})$}

In this section, we use the method of the weight space decomposition of $S$ with respect to the abelian subalgebra $T_S$ to prove that all skew-symmetric super-biderivations of $S$ are inner. For convenience, we use $A$ and $B$ denote $A(m,\underline{t})$ and $B(n)$.

Set $T_S=\Span_{\BF}\{D_{k_0k_1}(x_{k_0}x_{k_1})\mid k_0\in Y_0,k_1\in Y_1\}$. Obviously, $T_S$ is an abelian subalgebra of $S$. For any $D_{ij}(x^{(\alpha)}x^u)\in S$, we have
\[[D_{k_0k_1}(x_{k_0}x_{k_1}),D_{ij}(x^{(\alpha)}x^u)]=(\alpha_{k_0}+\delta_{(k_1\in  u)}-\delta_{ik_0}-\delta_{jk_0}-\delta_{ik_1}-\delta_{jk_1})D_{ij}(x^{(\alpha)}x^u),\]
where
\[\delta_{(\text{P})}=
\begin{cases}
	\ 1 \qquad \text{P is ture},\\
	\ 0 \qquad \text{P is false}.
\end{cases}\]
Fixed $\alpha\in A$, $u\in B$ and $i,j\in Y$, we define a linear function $(\alpha+\langle u\rangle+i+j):T_S\rightarrow \BF$
\[(\alpha+\langle u\rangle+i+j)(D_{k_0k_1}(x_{k_0}x_{k_1}))=\alpha_{k_0}+\delta_{(k_1\in  u)}-\delta_{ik_0}-\delta_{jk_0}-\delta_{ik_1}-\delta_{jk_1}.\]
Further, $S$ have a weight space decomposition with respect to $T_S$:
\[S=\bigoplus_{(\alpha+\langle u\rangle+i+j)}S_{(\alpha+\langle u\rangle+i+j)},\]
where
\begin{align*}
	S_{(\alpha+\langle u\rangle+i+j)}=\Span_{\BF}\{D_{st}(x^{(\beta)}x^v)\in S\mid &[D_{k_0k_1}(x_{k_0}x_{k_1}),D_{st}(x^{(\beta)}x^v)]=\\
	&(\alpha_{k_0}+\delta_{(k_1\in  u)}-\delta_{ik_0}-\delta_{jk_0}-\delta_{ik_1}-\delta_{jk_1})D_{st}(x^{(\beta)}x^v)\}.
\end{align*}

\begin{Lemma}\label{lem4.1}  	
	Let $\phi$ be a skew-symmetric super-biderivation of $S$. If $[x,y] = 0$ for $x,y \in S$, we have 
	\[\phi(x,y) = 0.\]	 
	   	  
	\begin{proof}
		Since $S$ is a simple Lie superalgebra, it is obvious that $S=[S,S]$ and $Z(S)=0$. By Lemma \ref{lem3.4}, if $[x,y]=0$ for $x, y \in S$, we obtain
		\[\phi(x,y)\in Z_S([S,S])=Z(S)=0.\]
		The proof is completed.
	\end{proof}
	
\end{Lemma}

\begin{Lemma}\label{lem4.2}		
	Let $\phi$ be a skew-symmetric super-biderivation of $S$. For any $D_{k_0k_1}(x_{k_0}x_{k_1})\in T_S$ and $D_{ij}(x^{(\alpha)}x^u)\in S$, we have
	\[\phi(D_{k_0k_1}(x_{k_0}x_{k_1}),D_{ij}(x^{(\alpha)}x^u))\in S_{(\alpha+\langle u\rangle+i+j)}.\]	
	
	\begin{proof}
		By lemma \ref{lem4.1}, it follows that $\phi(D_{k_0k_1}(x_{k_0}x_{k_1}),D_{l_0l_1}(x_{l_0}x_{l_1}))=0$ for any $k_0,l_0\in Y_0,k_1,l_1\in Y_1$ from $[D_{k_0k_1}(x_{k_0}x_{k_1}),D_{l_0l_1}(x_{l_0}x_{l_1})]=0$. Note that $\d(D_{k_0k_1}(x_{k_0}x_{k_1}))=\bar{0}$, then for any $D_{ij}(x^{(\alpha)}x^u)\in S$, it is clear that		
		\begin{align*}
			&[D_{l_0l_1}(x_{l_0}x_{l_1}),\phi(D_{k_0k_1}(x_{k_0}x_{k_1}),D_{ij}(x^{(\alpha)}x^u))]\\
			=&(-1)^{(\d(\phi)+\bar{0})\,\bar{0}}\big(\phi(D_{k_0k_1}(x_{k_0}x_{k_1}),[D_{l_0l_1}(x_{l_0}x_{l_1}),D_{ij}(x^{(\alpha)}x^u)])\\
			&-[\phi(D_{k_0k_1}(x_{k_0}x_{k_1}),D_{l_0l_1}(x_{l_0}x_{l_1})),D_{ij}(x^{(\alpha)}x^u)]\big)\\
			=& (\alpha_{l_0}+\delta_{(l_1\in  u)}-\delta_{il_0}-\delta_{jl_0}-\delta_{il_1}-\delta_{jl_1})\phi(D_{k_0k_1}(x_{k_0}x_{k_1}),D_{ij}(x^{(\alpha)}x^u)).
		\end{align*}
	    The proof is completed.
	\end{proof}
    
\end{Lemma}

We want to prove that all skew-symmetric super-biderivations of $S$ are inner. First, we need to prove the conclusion works on some elements of $S$. Therefore we give some specific weight spaces about the weight space decomposition of $S$ with respect to $T_S$.

\begin{Lemma}\label{lem4.3}
	Let $k_0,i_0,j_0\in Y_0,k_1,i_1,j_1\in Y_1$ and $i_0\neq k_0,j_0\neq i_0,i_1\neq k_1,j_1\neq i_1$. Then the following statements hold:
	\begin{align*}
		&(1)S_{(\varepsilon_{i_0}+k_0+i_0)}=\sum_{\alpha\in A,t\in Y_0\backslash \{k_0\}}\BF D_{k_0t}\big((\prod_{r\in Y_0\backslash \{t\}}x^{(\alpha_r^{\bar{0}}\varepsilon_r)})x^{(\alpha_t^{\bar{1}}\varepsilon_t)}\big)\\
		&\hspace*{8em}+\sum_{\alpha\in A,t\in Y_1}\BF D_{k_0t}\big((\prod_{r\in Y_0}x^{(\alpha_r^{\bar{0}}\varepsilon_r)})x_t\big);\\
		&(2)S_{(\langle i_1\rangle+k_1+i_1)}=\sum_{\alpha\in A,t\in Y_1}\BF D_{k_1t}\big((\prod_{r\in Y_0}x^{(\alpha_r^{\bar{0}}\varepsilon_r)})x_t\big);\\
		&(3)S_{(\varepsilon_{j_0}+\langle k_1\rangle+i_0+j_0)}=\sum_{\alpha\in A,t\in Y_0\backslash \{i_0\}}\BF D_{i_0t}\big((\prod_{r\in Y_0\backslash \{t\}}x^{(\alpha_r^{\bar{0}}\varepsilon_r)})x^{(\alpha_t^{\bar{1}}\varepsilon_t)}x_{k_1}\big)\\
		&\hspace*{10em}+\sum_{\alpha\in A,t\in Y_1\backslash\{k_1\}}\BF D_{i_0t}\big((\prod_{r\in Y_0}x^{(\alpha_r^{\bar{0}}\varepsilon_r)})x_tx_{k_1}\big);\\
		&(4)S_{(\varepsilon_{k_0}+\langle j_1\rangle+i_1+j_1)}=\sum_{\alpha\in A,t\in Y_1}\BF D_{i_1t}\big((\prod_{r\in Y_0\backslash \{k_0\}}x^{(\alpha_r^{\bar{0}}\varepsilon_r)})x^{(\alpha_{k_0}^{\bar{1}}\varepsilon_{k_0})}x_t\big);&
	\end{align*}
    where $\alpha_r^{\bar{q}}$ denotes some integer and $\alpha_r^{\bar{q}}\equiv q\pmod{p}$.
\end{Lemma}

\begin{Lemma}\label{lem4.4}
	\textup{\cite{zhang2005}}
	For any $f\in \Lambda(m,n;\underline{t})$ and $i,j,k\in Y$, we have
	\[[D_k,D_{ij}(f)]=(-1)^{\tau(k)\tau(i)}D_{ij}(D_k(f)).\]
\end{Lemma}

\begin{Lemma}\label{lem4.5}
	Let $\phi$ is a skew-symmetric super-biderivation of $S$. For any $k_0,i_0\in Y_0,k_1\in Y_1$ and $i_0\neq k_0$, there is an element $\lambda_{k_0k_1}\in\BF$ such that
	\[\phi(D_{k_0k_1}(x_{k_0}x_{k_1}),D_{k_0i_0}(x_{i_0}))=\lambda_{k_0k_1}[D_{k_0k_1}(x_{k_0}x_{k_1}),D_{k_0i_0}(x_{i_0})],\]
	where $\lambda_{k_0k_1}$ is dependent on $k_0$ and $k_1$.
	
	\begin{proof}
		By Lemma \ref{lem4.3} (1),  we can suppose that		
		
		\begin{equation}\label{eq4.1}
			\begin{aligned}
				 \phi(D_{k_0k_1}(x_{k_0}x_{k_1}),D_{k_0i_0}(x_{i_0}))
				&=\sum_{\alpha\in A,t\in Y_0\backslash \{k_0\}}a_0(\alpha,t,k_0,k_1)D_{k_0t}\big((\prod_{r\in Y_0\backslash \{t\}}x^{(\alpha_r^{\bar{0}}\varepsilon_r)})x^{(\alpha_t^{\bar{1}}\varepsilon_t)}\big)\\
				&\qquad+\sum_{\alpha\in A,t\in Y_1}a_1(\alpha,t,k_0,k_1)D_{k_0t}\big((\prod_{r\in Y_0}x^{(\alpha_r^{\bar{0}}\varepsilon_r)})x_t\big),
			\end{aligned}
		\end{equation}		
	    where $a_0(\alpha,t,k_0,k_1),a_1(\alpha,t,k_0,k_1)\in\BF$. By Lemma \ref{lem4.1}, \ref{lem4.4} and equation \eqref{eq4.1}, for any $l\in Y\backslash\{k_0,k_1\}$, we have
	    \begin{align*}
	    	0&=(-1)^{\d(\phi)\tau(l)}\big(\phi(D_{k_0k_1}(x_{k_0}x_{k_1}),[D_l,D_{k_0i_0}(x_{i_0})])-[\phi(D_{k_0k_1}(x_{k_0}x_{k_1}),D_l),D_{k_0i_0}(x_{i_0})]\big)\\
	    	&=[D_l,\phi(D_{k_0k_1}(x_{k_0}x_{k_1}),D_{k_0i_0}(x_{i_0}))]\\
	    	&=\sum_{\alpha\in A,t\in Y_0\backslash \{k_0\}}a_0(\alpha,t,k_0,k_1)(-1)^{\tau(l)\tau{(k_0)}}D_{k_0t}\big(D_l((\prod_{r\in Y_0\backslash \{t\}}x^{(\alpha_r^{\bar{0}}\varepsilon_r)})x^{(\alpha_t^{\bar{1}}\varepsilon_t)})\big)\\
	    	&\qquad+\sum_{\alpha\in A,t\in Y_1}a_1(\alpha,t,k_0,k_1)(-1)^{\tau(l)\tau{(k_0)}}D_{k_0t}\big(D_l((\prod_{r\in Y_0}x^{(\alpha_r^{\bar{0}}\varepsilon_r)})x_t)\big).
	    \end{align*}
        According to the above equation, we can draw the following conclusion. For the first summation formula, if $\alpha_r^{\bar{0}}>0$ for $r\neq k_0,t$ or $\alpha_t^{\bar{1}}>1$, we find that $a_0(\alpha,t,k_0,k_1)D_{k_0t}\big((\prod_{r\in Y_0\backslash \{t\}}x^{(\alpha_r^{\bar{0}}\varepsilon_r)})x^{(\alpha_t^{\bar{1}}\varepsilon_t)}\big)$ either has a coefficient of zero or can be eliminated from some terms in the second summation formula. In a word, this term will not exist in equation \eqref{eq4.1}. For the second summation formula, if $\alpha_r^{\bar{0}}>0$ for $r\neq k_0$, we find that $a_1(\alpha,t,k_0,k_1)D_{k_0t}\big((\prod_{r\in Y_0}x^{(\alpha_r^{\bar{0}}\varepsilon_r)})x_t\big)$ either has a coefficient of zero or can be eliminated from some terms in the first summation formula. Similarly, this term will not exist in equation \eqref{eq4.1}. Thus we can suppose that
        \begin{align*}
            \phi(D_{k_0k_1}(x_{k_0}x_{k_1}),D_{k_0i_0}(x_{i_0}))
        	&=\sum_{\alpha\in A,t\in Y_0\backslash \{k_0\}}a_0(\alpha,t,k_0,k_1)D_{k_0t}(x^{(\alpha_{k_0}^{\bar{0}}\varepsilon_{k_0})}x^{(\varepsilon_t)})\\
        	&\quad+\sum_{\alpha\in A,t\in Y_1}a_1(\alpha,t,k_0,k_1)D_{k_0t}(x^{(\alpha_{k_0}^{\bar{0}}\varepsilon_{k_0})}x_t).
        \end{align*}
        Since $\d(D_{k_0k_1}(x_{k_0}x_{k_1}))+\d(D_{k_0i_0}(x_{i_0}))=0$, by Lemma \ref{lem3.3}, we have
        \begin{align*}
        	0&=[\phi(D_{k_0k_1}(x_{k_0}x_{k_1}),D_{k_0i_0}(x_{i_0})),[D_{k_0k_1}(x_{k_0}x_{k_1}),D_{k_0i_0}(x_{i_0})]]\\
        	&=[\phi(D_{k_0k_1}(x_{k_0}x_{k_1}),D_{k_0i_0}(x_{i_0})),D_{k_0}]\\
        	&=-[D_{k_0},\phi(D_{k_0k_1}(x_{k_0}x_{k_1}),D_{k_0i_0}(x_{i_0}))].
        \end{align*}
        Similarly, by Lemma \ref{lem4.4}, it is easy to verify that
        \begin{align*}
        	&\quad\ \phi(D_{k_0k_1}(x_{k_0}x_{k_1}),D_{k_0i_0}(x_{i_0}))\\
        	&=\sum_{t\in Y_0\backslash \{k_0\}}a_0(t,k_0,k_1)D_{k_0t}(x^{(\varepsilon_t)})+\sum_{t\in Y_1}a_1(t,k_0,k_1)D_{k_0t}(x_t)\\
        	&=-\sum_{t\in Y_0\backslash \{k_0\}}a_0(t,k_0,k_1)D_{k_0}+\sum_{t\in Y_1}a_1(t,k_0,k_1)D_{k_0}\\
        	&:= a(k_0,k_1)D_{k_0}.
        \end{align*}
        Set $\lambda_{k_0k_1}=a(k_0,k_1)$. Since $[D_{k_0k_1}(x_{k_0}x_{k_1}),D_{k_0i_0}(x_{i_0})]=[x_{k_0}D_{k_0}+x_{k_1}D_{k_1},-D_{k_0}]=D_{k_0}$,
        we conclude that
        \[\phi(D_{k_0k_1}(x_{k_0}x_{k_1}),D_{k_0i_0}(x_{i_0}))=\lambda_{k_0k_1}[D_{k_0k_1}(x_{k_0}x_{k_1}),D_{k_0i_0}(x_{i_0})].\]
        where $\lambda_{k_0k_1}$ is dependent on $k_0$ and $k_1$.
	\end{proof}
    
\end{Lemma}

\begin{Lemma}
	All $\BZ_2$-homogeneous skew-symmetric super-biderivations of $S$ are even. 
	
	\begin{proof}
		Suppose $\phi$ is a $\BZ_2$-homogeneous skew-symmetric super-biderivation of $S$. Due to Lemma \ref{lem4.5}, we know that  $\phi(D_{k_0k_1}(x_{k_0}x_{k_1}),D_{k_0i_0}(x_{i_0}))$ and $[D_{k_0k_1}(x_{k_0}x_{k_1}),D_{k_0i_0}(x_{i_0})]$ have the same $\BZ_2$-degree, so $\phi$ and $[\cdot,\cdot]$ have the same $\BZ_2$-degree. Then $\phi$ is even.
	\end{proof}
	
\end{Lemma}

\begin{Corollary}
	All skew-symmetric super-biderivations of $S$ are $\BZ_2$-homogeneous, and fuether even. 
\end{Corollary}

\begin{Lemma}\label{lem4.8}
	Let $\phi$ is a skew-symmetric super-biderivation of $S$. For any $k_0\in Y_0,k_1,i_1\in Y_1$ and $i_1\neq k_1$, there is an element $\mu_{k_0k_1}\in\BF$ such that
	\[\phi(D_{k_0k_1}(x_{k_0}x_{k_1}),D_{k_1i_1}(x_{i_1}))=\mu_{k_0k_1}[D_{k_0k_1}(x_{k_0}x_{k_1}),D_{k_1i_1}(x_{i_1})],\]
	where $\mu_{k_0k_1}$ is dependent on $k_0$ and $k_1$.
	
	\begin{proof}
		By Lemma \ref{lem4.3} (2),  we can suppose that		
			\begin{align}\label{eq4.2}
				\phi(D_{k_0k_1}(x_{k_0}x_{k_1}),D_{k_1i_1}(x_{i_1}))
				=\sum_{\alpha\in A,t\in Y_1}a_1(\alpha,t,k_0,k_1)D_{k_1t}\big((\prod_{r\in Y_0}x^{(\alpha_r^{\bar{0}}\varepsilon_r)})x_t\big),
			\end{align}			
		where $a_1(\alpha,t,k_0,k_1)\in\BF$. By Lemma \ref{lem4.1}, \ref{lem4.4} and equation \eqref{eq4.2}, for any $l\in Y\backslash\{k_0,k_1\}$, we have
		\begin{align*}
			0&=\phi(D_{k_0k_1}(x_{k_0}x_{k_1}),[D_l,D_{k_1i_1}(x_{i_1})])-[\phi(D_{k_0k_1}(x_{k_0}x_{k_1}),D_l),D_{k_1i_1}(x_{i_1})]\\
			&=[D_l,\phi(D_{k_0k_1}(x_{k_0}x_{k_1}),D_{k_1i_1}(x_{i_1}))]\\
			&=\sum_{\alpha\in A,t\in Y_1}a_1(\alpha,t,k_0,k_1)(-1)^{\tau(l)\tau{(k_1)}}D_{k_1t}\big(D_l((\prod_{r\in Y_0}x^{(\alpha_r^{\bar{0}}\varepsilon_r)})x_t)\big).
		\end{align*}
		According to the above equation, we can draw the following conclusion. If $\alpha_r^{\bar{0}}>0$ for $r\neq k_0$, we find that $\sum_{t\in Y_1}(\delta_{k_1t}+1)a_1(\alpha,t,k_0,k_1)=0$. This is equivalent to the fact that $a_1(\alpha,t,k_0,k_1)D_{k_1t}\big((\prod_{r\in Y_0}x^{(\alpha_r^{\bar{0}}\varepsilon_r)})x_t\big)$ will be eliminated in equation \eqref{eq4.2}. Thus we can suppose that
		\begin{align*}
			\phi(D_{k_0k_1}(x_{k_0}x_{k_1}),D_{k_1i_1}(x_{i_1}))
			=\sum_{\alpha\in A,t\in Y_1}a_1(\alpha,t,k_0,k_1)D_{k_1t}(x^{(\alpha_{k_0}^{\bar{0}}\varepsilon_{k_0})}x_t).			
		\end{align*}
		By Lemma \ref{lem3.2} and \ref{lem4.5}, we have
		\begin{align*}
			0&=[\phi(D_{k_0k_1}(x_{k_0}x_{k_1}),D_{k_0i_0}(x_{i_0})),[D_{k_0k_1}(x_{k_0}x_{k_1}),D_{k_1i_1}(x_{i_1})]]\\
			&=[[D_{k_0k_1}(x_{k_0}x_{k_1}),D_{k_0i_0}(x_{i_0})],\phi(D_{k_0k_1}(x_{k_0}x_{k_1}),D_{k_1i_1}(x_{i_1}))]\\
			&=[D_{k_0},\phi(D_{k_0k_1}(x_{k_0}x_{k_1}),D_{k_1i_1}(x_{i_1}))].
		\end{align*}
		Similarly, by Lemma \ref{lem4.4}, it is easy to verify that
		\begin{align*}
			\phi(D_{k_0k_1}(x_{k_0}x_{k_1}),D_{k_1i_1}(x_{i_1}))
			&=\sum_{t\in Y_1}a_1(t,k_0,k_1)D_{k_1t}(x_t)\\
			&=-\sum_{t\in Y_1\backslash\{k_1\}}a_1(t,k_0,k_1)D_{k_1}-2a_1(k_1,k_0,k_1)D_{k_1}\\
			&:= a(k_0,k_1)D_{k_1}.
		\end{align*}
		Set $\mu_{k_0k_1}=a(k_0,k_1)$. Since $[D_{k_0k_1}(x_{k_0}x_{k_1}),D_{k_1i_1}(x_{i_1})]=[x_{k_0}D_{k_0}+x_{k_1}D_{k_1},-D_{k_1}]=D_{k_1}$,
		we conclude that
		\[\phi(D_{k_0k_1}(x_{k_0}x_{k_1}),D_{k_1i_1}(x_{i_1}))=\mu_{k_0k_1}[D_{k_0k_1}(x_{k_0}x_{k_1}),D_{k_1i_1}(x_{i_1})].\]
		where $\mu_{k_0k_1}$ is dependent on $k_0$ and $k_1$.
	\end{proof}
	
\end{Lemma}

\begin{Lemma}\label{lem4.9}
	Let $\phi$ is a skew-symmetric super-biderivation of $S$. For any $k_0,i_0,j_0\in Y_0,k_1\in Y_1$ and $i_0\neq k_0,j_0\neq i_0$, there is an element $\lambda_{k_0k_1i_0}\in\BF$ such that
	\[\phi(D_{k_0k_1}(x_{k_0}x_{k_1}),D_{i_0j_0}(x_{j_0}x_{k_1}))=\lambda_{k_0k_1i_0}[D_{k_0k_1}(x_{k_0}x_{k_1}),D_{i_0j_0}(x_{j_0}x_{k_1})],\]
	where $\lambda_{k_0k_1i_0}$ is dependent on $k_0$, $k_1$ and $i_0$.
	
	\begin{proof}
		By Lemma \ref{lem4.3} (3),  we can suppose that		
		\begin{equation}\label{eq4.3}
			\begin{aligned}
				&\quad\ \phi(D_{k_0k_1}(x_{k_0}x_{k_1}),D_{i_0j_0}(x_{j_0}x_{k_1}))\\
				&=\sum_{\alpha\in A,t\in Y_0\backslash \{i_0\}}a_0(\alpha,t,k_0,k_1,i_0)D_{i_0t}\big((\prod_{r\in Y_0\backslash \{t\}}x^{(\alpha_r^{\bar{0}}\varepsilon_r)})x^{(\alpha_t^{\bar{1}}\varepsilon_t)}x_{k_1}\big)\\
				&\qquad+\sum_{\alpha\in A,t\in Y_1\backslash\{k_1\}}a_1(\alpha,t,k_0,k_1,i_0)D_{i_0t}\big((\prod_{r\in Y_0}x^{(\alpha_r^{\bar{0}}\varepsilon_r)})x_tx_{k_1}\big),
			\end{aligned}
		\end{equation}		
		where $a_0(\alpha,t,k_0,k_1,i_0),a_1(\alpha,t,k_0,k_1,i_0)\in\BF$. By Lemma \ref{lem4.1}, \ref{lem4.4} and equation \eqref{eq4.3}, for any $l\in Y\backslash\{k_0,k_1\}$, we have
		\begin{align*}
			0&=\phi(D_{k_0k_1}(x_{k_0}x_{k_1}),[D_l,D_{i_0j_0}(x_{j_0}x_{k_1})])-[\phi(D_{k_0k_1}(x_{k_0}x_{k_1}),D_l),D_{i_0j_0}(x_{j_0}x_{k_1})]\\
			&=[D_l,\phi(D_{k_0k_1}(x_{k_0}x_{k_1}),D_{i_0j_0}(x_{j_0}x_{k_1}))]\\
			&=\sum_{\alpha\in A,t\in Y_0\backslash \{i_0\}}a_0(\alpha,t,k_0,k_1,i_0)(-1)^{\tau(l)\tau{(i_0)}}D_{i_0t}\big(D_l((\prod_{r\in Y_0\backslash \{t\}}x^{(\alpha_r^{\bar{0}}\varepsilon_r)})x^{(\alpha_t^{\bar{1}}\varepsilon_t)}x_{k_1})\big)\\
			&\qquad+\sum_{\alpha\in A,t\in Y_1\backslash\{k_1\}}a_1(\alpha,t,k_0,k_1,i_0)(-1)^{\tau(l)\tau{(i_0)}}D_{i_0t}\big(D_l((\prod_{r\in Y_0}x^{(\alpha_r^{\bar{0}}\varepsilon_r)})x_tx_{k_1})\big).
		\end{align*}
		According to the above equation, we can draw the following conclusion. For the first summation formula, if $\alpha_r^{\bar{0}}>0$ for $r\neq k_0,t$ or $\alpha_t^{\bar{1}}>1$, we find that $a_0(\alpha,t,k_0,k_1,i_0)D_{i_0t}\big((\prod_{r\in Y_0\backslash \{t\}}x^{(\alpha_r^{\bar{0}}\varepsilon_r)})x^{(\alpha_t^{\bar{1}}\varepsilon_t)}x_{k_1}\big)$ either has a coefficient of zero or can be eliminated from some terms in the second summation formula. In a word, this term will not exist in equation \eqref{eq4.3}. For the second summation formula, if $\alpha_r^{\bar{0}}>0$ for $r\neq k_0$, we find that $a_1(\alpha,t,k_0,k_1,i_0)D_{i_0t}\big((\prod_{r\in Y_0}x^{(\alpha_r^{\bar{0}}\varepsilon_r)})x_tx_{k_1}\big)$ either has a coefficient of zero or can be eliminated from some terms in the first summation formula. Similarly, this term will not exist in equation \eqref{eq4.3}. Thus we can suppose that
		\begin{align*}
			 \phi(D_{k_0k_1}(x_{k_0}x_{k_1}),D_{i_0j_0}(x_{j_0}x_{k_1}))
			&=\sum_{\alpha\in A,t\in Y_0\backslash \{i_0\}}a_0(\alpha,t,k_0,k_1,i_0)D_{i_0t}(x^{(\alpha_{k_0}^{\bar{0}}\varepsilon_{k_0})}x^{(\varepsilon_t)}x_{k_1})\\
			&\quad+\sum_{\alpha\in A,t\in Y_1\backslash\{k_1\}}a_1(\alpha,t,k_0,k_1,i_0)D_{i_0t}(x^{(\alpha_{k_0}^{\bar{0}}\varepsilon_{k_0})}x_tx_{k_1}).
		\end{align*}
		By Lemma \ref{lem3.2} and \ref{lem4.5}, we have
		\begin{align*}
			0&=[\phi(D_{k_0k_1}(x_{k_0}x_{k_1}),D_{k_0i_0}(x_{i_0})),[D_{k_0k_1}(x_{k_0}x_{k_1}),D_{i_0j_0}(x_{j_0}x_{k_1})]]\\
			&=[[D_{k_0k_1}(x_{k_0}x_{k_1}),D_{k_0i_0}(x_{i_0})],\phi(D_{k_0k_1}(x_{k_0}x_{k_1}),D_{i_0j_0}(x_{j_0}x_{k_1}))]\\
			&=[D_{k_0},\phi(D_{k_0k_1}(x_{k_0}x_{k_1}),D_{i_0j_0}(x_{j_0}x_{k_1}))].
		\end{align*}
		Similarly, by Lemma \ref{lem4.4}, it is easy to verify that
		\begin{align*}
			&\quad\ \phi(D_{k_0k_1}(x_{k_0}x_{k_1}),D_{i_0j_0}(x_{j_0}x_{k_1}))\\
			&=\sum_{t\in Y_0\backslash \{i_0\}}a_0(t,k_0,k_1,i_0)D_{i_0t}(x^{(\varepsilon_t)}x_{k_1})+\sum_{t\in Y_1\backslash\{k_1\}}a_1(t,k_0,k_1,i_0)D_{i_0t}(x_tx_{k_1})\\
			&=-\sum_{t\in Y_0\backslash \{i_0\}}a_0(t,k_0,k_1,i_0)x_{k_1}D_{i_0}-\sum_{t\in Y_1\backslash\{k_1\}}a_1(t,k_0,k_1,i_0)x_{k_1}D_{i_0}\\
			&:= a(k_0,k_1,i_0)x_{k_1}D_{i_0}.
		\end{align*}
		Set $\lambda_{k_0k_1i_0}=-a(k_0,k_1,i_0)$. Since $$[D_{k_0k_1}(x_{k_0}x_{k_1}),D_{i_0j_0}(x_{j_0}x_{k_1})]=[x_{k_0}D_{k_0}+x_{k_1}D_{k_1},-x_{k_1}D_{i_0}]=-x_{k_1}D_{i_0},$$
		we conclude that
		\[\phi(D_{k_0k_1}(x_{k_0}x_{k_1}),D_{i_0j_0}(x_{j_0}x_{k_1}))=\lambda_{k_0k_1i_0}[D_{k_0k_1}(x_{k_0}x_{k_1}),D_{i_0j_0}(x_{j_0}x_{k_1})].\]
		where $\lambda_{k_0k_1i_0}$ is dependent on $k_0$, $k_1$ and $i_0$.
	\end{proof}
	
\end{Lemma}

\begin{Lemma}\label{lem4.10}
	Let $\phi$ is a skew-symmetric super-biderivation of $S$. For any $k_0,\in Y_0,k_1,i_1,j_1\in Y_1$ and $i_1\neq k_1,j_1\neq i_1$, there is an element $\mu_{k_0k_1i_1}\in\BF$ such that
	\[\phi(D_{k_0k_1}(x_{k_0}x_{k_1}),D_{i_1j_1}(x_{j_1}x_{k_0}))=\mu_{k_0k_1i_1}[D_{k_0k_1}(x_{k_0}x_{k_1}),D_{i_1j_1}(x_{j_1}x_{k_0})],\]
	where $\mu_{k_0k_1i_1}$ is dependent on $k_0$, $k_1$ and $i_1$.
	
	\begin{proof}
		By Lemma \ref{lem4.3} (4),  we can suppose that		
		\begin{equation}\label{eq4.4}
			\begin{aligned}
				&\quad\ \phi(D_{k_0k_1}(x_{k_0}x_{k_1}),D_{i_1j_1}(x_{j_1}x_{k_0}))\\
				&=\sum_{\alpha\in A,t\in Y_1}a_1(\alpha,t,k_0,k_1,i_1)D_{i_1t}\big((\prod_{r\in Y_0\backslash \{k_0\}}x^{(\alpha_r^{\bar{0}}\varepsilon_r)})x^{(\alpha_{k_0}^{\bar{1}}\varepsilon_{k_0})}x_t\big),
			\end{aligned}
		\end{equation}		
		where $a_1(\alpha,t,k_0,k_1,i_1)\in\BF$. By Lemma \ref{lem4.1}, \ref{lem4.4} and equation \eqref{eq4.4}, for any $l\in Y\backslash\{k_0,k_1\}$, we have
		\begin{align*}
			0&=\phi(D_{k_0k_1}(x_{k_0}x_{k_1}),[D_l,D_{i_1j_1}(x_{j_1}x_{k_0})])-[\phi(D_{k_0k_1}(x_{k_0}x_{k_1}),D_l),D_{i_1j_1}(x_{j_1}x_{k_0})]\\
			&=[D_l,\phi(D_{k_0k_1}(x_{k_0}x_{k_1}),D_{i_1j_1}(x_{j_1}x_{k_0}))]\\
			&=\sum_{\alpha\in A,t\in Y_1}a_1(\alpha,t,k_0,k_1,i_1)(-1)^{\tau(l)\tau{(i_1)}}D_{i_1t}\big(D_l((\prod_{r\in Y_0\backslash \{k_0\}}x^{(\alpha_r^{\bar{0}}\varepsilon_r)})x^{(\alpha_{k_0}^{\bar{1}}\varepsilon_{k_0})}x_t)\big).
		\end{align*}
		According to the above equation, we can draw the following conclusion. If $\alpha_r^{\bar{0}}>0$ for $r\neq k_0$, we find that $\sum_{t\in Y_1}(\delta_{i_1t}+1)a_1(\alpha,t,k_0,k_1,i_1)=0$. This is equivalent to the fact that $a_1(\alpha,t,k_0,k_1,i_1)D_{i_1t}\big((\prod_{r\in Y_0\backslash \{k_0\}}\\x^{(\alpha_r^{\bar{0}}\varepsilon_r)})x^{(\alpha_{k_0}^{\bar{1}}\varepsilon_{k_0})}x_t\big)$ will be eliminated in equation \eqref{eq4.4}. Thus we can suppose that
		\begin{align*}
			\phi(D_{k_0k_1}(x_{k_0}x_{k_1}),D_{i_1j_1}(x_{j_1}x_{k_0}))
			=\sum_{t\in Y_1}a_1(\alpha,t,k_0,k_1,i_1)D_{i_1t}(x^{(\alpha_{k_0}^{\bar{1}}\varepsilon_{k_0})}x_t).
		\end{align*}
		By Lemma \ref{lem3.2} and \ref{lem4.5}, we have
		\begin{align*}
			-\lambda_{k_0k_1}D_{i_1}&=[\lambda_{k_0k_1}D_{k_0},-x_{k_0}D_{i_1}]\\
			&=[\phi(D_{k_0k_1}(x_{k_0}x_{k_1}),D_{k_0i_0}(x_{i_0})),[D_{k_0k_1}(x_{k_0}x_{k_1}),D_{i_1j_1}(x_{j_1}x_{k_0})]]\\
			&=[[D_{k_0k_1}(x_{k_0}x_{k_1}),D_{k_0i_0}(x_{i_0})],\phi(D_{k_0k_1}(x_{k_0}x_{k_1}),D_{i_1j_1}(x_{j_1}x_{k_0}))]\\
			&=[D_{k_0},\phi(D_{k_0k_1}(x_{k_0}x_{k_1}),D_{i_1j_1}(x_{j_1}x_{k_0}))].
		\end{align*}
	    By Lemma \ref{lem4.4}, we have
	    \begin{align*}
	    	-\lambda_{k_0k_1}D_{i_1}
	    	&=\sum_{t\in Y_1}a_1(\alpha,t,k_0,k_1,i_1)D_{i_1t}\big(D_{k_0}(x^{(\alpha_{k_0}^{\bar{1}}\varepsilon_{k_0})}x_t)\big)\\
	    	&=\sum_{t\in Y_1}a_1(\alpha,t,k_0,k_1,i_1)D_{i_1t}(x^{(\alpha_{k_0}^{\bar{0}}\varepsilon_{k_0})}x_t)\\
	    	&=-\big(\sum_{t\in Y_1\backslash\{i_1\}}a_1(\alpha,t,k_0,k_1,i_1)+2a_1(\alpha,i_1,k_0,k_1,i_1)\big)x^{(\alpha_{k_0}^{\bar{0}}\varepsilon_{k_0})}D_{i_1}.
	    \end{align*}		
		By computing the equation,  it is easy to verify that
		\begin{align*}
			\phi(D_{k_0k_1}(x_{k_0}x_{k_1}),D_{i_1j_1}(x_{j_1}x_{k_0}))
			&=\sum_{t\in Y_1}a_1(t,k_0,k_1,i_1)D_{i_1t}(x^{(\varepsilon_{k_0})}x_t)\\
			&=-\big(\sum_{t\in Y_1\backslash \{i_1\}}a_1(t,k_0,k_1,i_1)+2a_1(i_1,k_0,k_1,i_1)\big)x_{k_0}D_{i_1}\\
			&:= a(k_0,k_1,i_1)x_{k_0}D_{i_1}.
		\end{align*}
		Set $\mu_{k_0k_1i_1}=-a(k_0,k_1,i_1)$. Since $$[D_{k_0k_1}(x_{k_0}x_{k_1}),D_{i_1j_1}(x_{j_1}x_{k_0})]=[x_{k_0}D_{k_0}+x_{k_1}D_{k_1},-x_{k_0}D_{i_1}]=-x_{k_0}D_{i_1},$$
		we conclude that
		\[\phi(D_{k_0k_1}(x_{k_0}x_{k_1}),D_{i_1j_1}(x_{j_1}x_{k_0}))=\mu_{k_0k_1i_1}[D_{k_0k_1}(x_{k_0}x_{k_1}),D_{i_1j_1}(x_{j_1}x_{k_0})].\]
		where $\mu_{k_0k_1i_1}$ is dependent on $k_0$, $k_1$ and $i_1$.
	\end{proof}
	
\end{Lemma}

\begin{Remark}
	The coefficients $\lambda_{k_0k_1}$, $\mu_{k_0k_1}$, $\lambda_{k_0k_1i_0}$ and $\mu_{k_0k_1i_1}$ in Lemmas \ref{lem4.5} $\sim$ \ref{lem4.10} are the same. This coefficient is dependent only on $\phi$, independent of $k_0$, $k_1$, $i_0$, $i_1$.
	
	\begin{proof}
		By Lemma \ref{lem3.2}, \ref{lem4.5} and \ref{lem4.10}, we have 
		\begin{align*}
			&\quad\ [\lambda_{k_0k_1}D_{k_0},-x_{k_0}D_{i_1}]\\			
			&=[\phi(D_{k_0k_1}(x_{k_0}x_{k_1}),D_{k_0i_0}(x_{i_0})),[D_{k_0k_1}(x_{k_0}x_{k_1}),D_{i_1j_1}(x_{j_1}x_{k_0})]]\\
			&=[[D_{k_0k_1}(x_{k_0}x_{k_1}),D_{k_0i_0}(x_{i_0})],\phi(D_{k_0k_1}(x_{k_0}x_{k_1}),D_{i_1j_1}(x_{j_1}x_{k_0}))]\\
			&=[D_{k_0},-\mu_{k_0k_1i_1}x_{k_0}D_{i_1}].
		\end{align*}
	    By direct calculation, it is easily seen that $(\lambda_{k_0k_1}-\mu_{k_0k_1i_1})D_{i_1}=0$. Since $D_{i_1}\neq 0$, we have that $\mu_{k_0k_1i_1}=\lambda_{k_0k_1}$. Similarly, by Lemma \ref{lem3.2}, \ref{lem4.8} and \ref{lem4.9}, we have that $\lambda_{k_0k_1i_0}=\mu_{k_0k_1}$. Then the conclusions of Lemma \ref{lem4.5} $\sim$ \ref{lem4.10} can be rewritten as
	    \begin{align}
	    	\phi(D_{k_0k_1}(x_{k_0}x_{k_1}),D_{k_0i_0}(x_{i_0}))&=\lambda_{k_0k_1}[D_{k_0k_1}(x_{k_0}x_{k_1}),D_{k_0i_0}(x_{i_0})],\label{eq4.5}\\
	    	\phi(D_{k_0k_1}(x_{k_0}x_{k_1}),D_{k_1i_1}(x_{i_1}))&=\mu_{k_0k_1}[D_{k_0k_1}(x_{k_0}x_{k_1}),D_{k_1i_1}(x_{i_1})],\label{eq4.6}\\
	    	\phi(D_{k_0k_1}(x_{k_0}x_{k_1}),D_{i_0j_0}(x_{j_0}x_{k_1}))&=\mu_{k_0k_1}[D_{k_0k_1}(x_{k_0}x_{k_1}),D_{i_0j_0}(x_{j_0}x_{k_1})],\label{eq4.7}\\
	    	\phi(D_{k_0k_1}(x_{k_0}x_{k_1}),D_{i_1j_1}(x_{j_1}x_{k_0}))&=\lambda_{k_0k_1}[D_{k_0k_1}(x_{k_0}x_{k_1}),D_{i_1j_1}(x_{j_1}x_{k_0})].\label{eq4.8}
	    \end{align}
        The equation \eqref{eq4.8} is equivalent to that 
        \begin{align}
        	\phi(D_{i_0i_1}(x_{i_0}x_{i_1}),D_{k_1j_1}(x_{j_1}x_{i_0}))=\lambda_{i_0i_1}[D_{i_0i_1}(x_{i_0}x_{i_1}),D_{k_1j_1}(x_{j_1}x_{i_0})].\label{eq4.9}
        \end{align}
        where $\lambda_{i_0i_1}$ is dependent on $i_0$ and $i_1$. By Lemma \ref{lem3.2}, equation \eqref{eq4.7} and \eqref{eq4.9}, we have 
        \begin{align*}
        	&\quad\ [-\mu_{k_0k_1}x_{k_1}D_{i_0},-x_{i_0}D_{k_1}]\\			
        	&=[\phi(D_{k_0k_1}(x_{k_0}x_{k_1}),D_{i_0j_0}(x_{j_0}x_{k_1})),[D_{i_0i_1}(x_{i_0}x_{i_1}),D_{k_1j_1}(x_{j_1}x_{i_0})]]\\
        	&=[[D_{k_0k_1}(x_{k_0}x_{k_1}),D_{i_0j_0}(x_{j_0}x_{k_1})],\phi(D_{i_0i_1}(x_{i_0}x_{i_1}),D_{k_1j_1}(x_{j_1}x_{i_0}))]\\
        	&=[-x_{k_1}D_{i_0},-\lambda_{i_0i_1}x_{i_0}D_{k_1}].
        \end{align*}
        By direct calculation, it is easily seen that $(\mu_{k_0k_1}-\lambda_{i_0i_1})(x_{k_1}D_{k_1}+x_{i_0}D_{i_0})=0$. Since $x_{k_1}D_{k_1}+x_{i_0}D_{i_0}\neq 0$, we have that $\mu_{k_0k_1}=\lambda_{i_0i_1}$. Note that $k_0$, $i_0$ are any two disparate elements in $Y_0$, and $k_1$, $i_1$ are any two disparate elements in $Y_1$. Then we can conclude that the coefficients in Lemmas \ref{lem4.5} $\sim$ \ref{lem4.10} are the same and dependent only on $\phi$.	    
	\end{proof}

\end{Remark}

Without loss of generality, we set $\lambda_{k_0k_1}=\mu_{k_0k_1}=\lambda_{k_0k_1i_0}=\mu_{k_0k_1i_1}=\nu$. Then the conclusions of Lemma \ref{lem4.5} $\sim$ \ref{lem4.10} can be rewritten as
    \begin{align}
    	\phi(D_{k_0k_1}(x_{k_0}x_{k_1}),D_{k_0i_0}(x_{i_0}))&=\nu[D_{k_0k_1}(x_{k_0}x_{k_1}),D_{k_0i_0}(x_{i_0})],\label{eq4.10}\\
    	\phi(D_{k_0k_1}(x_{k_0}x_{k_1}),D_{k_1i_1}(x_{i_1}))&=\nu[D_{k_0k_1}(x_{k_0}x_{k_1}),D_{k_1i_1}(x_{i_1})],\label{eq4.11}\\
    	\phi(D_{k_0k_1}(x_{k_0}x_{k_1}),D_{i_0j_0}(x_{j_0}x_{k_1}))&=\nu[D_{k_0k_1}(x_{k_0}x_{k_1}),D_{i_0j_0}(x_{j_0}x_{k_1})],\label{eq4.12}\\
    	\phi(D_{k_0k_1}(x_{k_0}x_{k_1}),D_{i_1j_1}(x_{j_1}x_{k_0}))&=\nu[D_{k_0k_1}(x_{k_0}x_{k_1}),D_{i_1j_1}(x_{j_1}x_{k_0})],\label{eq4.13}
    \end{align}
where $k_0,i_0,j_0\in Y_0$, $k_1,i_1,j_1\in Y_1$ and $i_0\neq k_0,j_0\neq i_0$, $i_1\neq k_1,j_1\neq i_1$.

\begin{Theorem}
	Let $S$ be the special Lie superalgebra $S(m,n;\underline{t})$ over the basic field $\BF$ of characteristic $p>2$. Then all skew-symmetric super-biderivations of $S$ are inner, i.e.
	\[\BDer(S)=\IBDer(S).\]
	
	\begin{proof}
		Suppose $\phi$ is a skew-symmetric super-biderivation of $S$. For any $D_{ij}(f),D_{st}(g)\in S$, by Lemma \ref{lem3.2} and equation \eqref{eq4.10}, we have
		\begin{align*}
			&\quad\ [\nu D_{k_0},[D_{ij}(f),D_{st}(g)]]\\			
			&=[\phi(D_{k_0k_1}(x_{k_0}x_{k_1}),D_{k_0i_0}(x_{i_0})),[D_{ij}(f),D_{st}(g)]]\\
			&=[[D_{k_0k_1}(x_{k_0}x_{k_1}),D_{k_0i_0}(x_{i_0})],\phi(D_{ij}(f),D_{st}(g))]\\
			&=[D_{k_0},\phi(D_{ij}(f),D_{st}(g))].
		\end{align*}
	    By direct calculation, it is easily seen that
	        \[0=[D_{k_0},\phi(D_{ij}(f),D_{st}(g))-\nu[D_{ij}(f),D_{st}(g)]].\]
	    Similarly, by Lemma \ref{lem3.2} and equation \eqref{eq4.11}, we have
	        \[0=[D_{k_1},\phi(D_{ij}(f),D_{st}(g))-\nu[D_{ij}(f),D_{st}(g)]].\]
	    Since $Z_{S_{-1}}(S)=\{E\in S\mid [E,D_i]=0,\forall\ i\in Y\}=S_{-1}$, we have that
	        \[\phi(D_{ij}(f),D_{st}(g))-\nu[D_{ij}(f),D_{st}(g)]=\sum_{l\in Y}b_lD_l.\]
	    By Lemma \ref{lem3.2} and equation \eqref{eq4.12}, we have
	    \begin{align*}
	    	&\quad\ [-\nu x_{k_1}D_{i_0},[D_{ij}(f),D_{st}(g)]]\\			
	    	&=[\phi(D_{k_0k_1}(x_{k_0}x_{k_1}),D_{i_0j_0}(x_{j_0}x_{k_1})),[D_{ij}(f),D_{st}(g)]]\\
	    	&=[[D_{k_0k_1}(x_{k_0}x_{k_1}),D_{i_0j_0}(x_{j_0}x_{k_1})],\phi(D_{ij}(f),D_{st}(g))]\\
	    	&=[-x_{k_1}D_{i_0},\phi(D_{ij}(f),D_{st}(g))].
	    \end{align*}
         By direct calculation, it is easily seen that
             \[0=[x_{k_1}D_{i_0},\sum_{l\in Y}b_lD_l]=b_{k_1}D_{i_0}.\]
         Since $D_{i_0}\neq 0$, we have that $b_{k_1}=0$. Similarly, by Lemma \ref{lem3.2} and equation \eqref{eq4.13}, we have
             \[0=[x_{k_0}D_{i_1},\sum_{l\in Y}b_lD_l]=-b_{k_0}D_{i_1}.\] 
         Since $D_{i_1}\neq 0$, we have that $b_{k_0}=0$. Then $b_l=0,\ \forall\ l\in Y$. Hence, for any $D_{ij}(f),D_{st}(g)\in S$, we have
             \[\phi(D_{ij}(f),D_{st}(g))=\nu[D_{ij}(f),D_{st}(g)].\]
         Thus $\phi$ is an inner super-biderivation. From the arbitrariness of $\phi$, the proof is complete.
	\end{proof}

\end{Theorem}

%

\end{document}